\documentclass[preprint,11pt]{article}
\usepackage{amsmath,amssymb}
\usepackage{multirow}
\usepackage{enumerate}
\usepackage{array}
\usepackage{enumerate}
\usepackage{algorithm}
\usepackage{algpseudocode}
\usepackage{algorithmicx}
\usepackage{graphicx}

\newtheorem{theo}{Theorem}[section]

\newtheorem{coro}[theo]{Corollary}
\newtheorem{lemm}[theo]{Lemma}
\newtheorem{prop}[theo]{Proposition}

\newtheorem{ejemplo}[theo]{Example}

\usepackage{xcolor}

\begin{document}

\begin{center}
{\Large \bf Classification of asexual diploid organisms by means of strongly isotopic evolution algebras defined over any field}
\end{center}

\begin{center}
{\large \em O. J. Falc\'on$^1$,\ R. M. Falc\'on$^2$\, J. N\'u\~nez$^{3}$}
\vspace{0.5cm}

{\small $^{1, 3}$ Faculty of Mathematics, University of Seville, Spain.\\
$^2$School of Building Engineering, University of Seville, Spain.\\
E-mail: {\em $^1$oscfalgan@yahoo.es,\ $^2$rafalgan@us.es,\ $^3$jnvaldes@us.es}}
\end{center}

\vspace{0.5cm}

\noindent {\large \bf Abstract.} Evolution algebras were introduced into Genetics to deal with the mechanism of inheritance of asexual organisms. Their distribution into isotopism classes is uniquely related with the mutation of alleles in non-Mendelian Genetics. This paper deals with such a distribution by means of Computational Algebraic Geometry. We focus in particular on the two-dimensional case, which is related to the asexual reproduction processes of diploid organisms. Specifically, we determine the existence of four isotopism classes, whatever the base field is, and we characterize the corresponding isomorphism classes.

\vspace{0.5cm}


\section{Introduction}

In the middle of the twentieth century, nonassociative algebras were introduced in Genetics by Etherington \cite{Etherington1939, Etherington1941, Etherington1941a} in order to endow Mendel's laws with a mathematical formulation that simulates the sexual reproduction and the mechanism of inheritance of an organism by considering the fusion of gametes into a zygote as an algebraic multiplication whose structure constants determine the probability distribution of the gametic output. Much more recently, in order to deal with asexual reproduction processes, Tian and Vojtechovsky \cite{Tian2008, Tian2006} introduced evolution algebras as a type of genetic algebra that makes possible to deal algebraically with the self-reproduction of alleles in non-Mendelian Genetics. The fundamentals of such algebras have been being developed in the last years with no probabilistic restrictions on the structure constants \cite{Cabrera2016, Camacho2013a, Casas2011, Casas2013, Dzhumadildaev2014, Dzhumadildaev2016, Labra2014, Ladra2013, Ladra2014, Omirov2015}. Nowadays, evolution algebras also constitute a fundamental connection between algebra, dynamic systems, Markov processes, Knot Theory, Graph Theory and Group Theory \cite{Khudoyberdiyev2015, Tian2008}.

\vspace{0.1cm}

A main problem in the theory of evolution algebras is their distribution into isomorphism and isotopism classes. On the one hand, the mentioned distribution into isomorphism classes has already been dealt with for two-dimensional evolution algebras over the complex field \cite{Camacho2013, Casas2014} and for nilpotent evolution algebras of dimension up to four over arbitrary fields \cite{Hegazi2015}. On the other hand, isotopisms have emerged as an interesting tool to simulate mutations in genetic algebras. In this regard, Holgate and Campos \cite{Campos1987, Holgate1966} had already considered isotopisms of genetic algebras as a way to formulate mathematically the mutation of alleles in the inheritance process. They showed indeed that certain known families of genetic algebras are isotopic. Nevertheless, to the best of the authors knowledge, isotopisms have not yet been considered in case of dealing with evolution algebras. The main goal of this paper is to delve further into this aspect.

\vspace{0.1cm}

The paper is organized as follows. In Section 2, we indicate some preliminaries concepts and results on isotopisms of algebras, genetic algebras and Computational Algebraic Geometry which will be used in the rest of the paper. Section 3 deals with the distribution of finite-dimensional evolution algebras over any base field into isotopism classes according to their structure tuples and to the dimension of their annihilators. Particularly, we determine the existence of four isotopism classes of two-dimensional evolution algebras, whatever the base field is. After that, we focus in Section 4 on the corresponding distribution of two-dimensional evolution algebras over any base field into isomorphism classes.

\section{Preliminaries}

In this section we expose some basic concepts and results on isotopisms of algebras, genetic algebras and Computational Algebraic Geometry that are used throughout the paper. For more details about these topics we refer, respectively, to the manuscripts of Albert \cite{Albert1942}, W\"orz-Busekros \cite{Worz1980}, Tian \cite{Tian2008} and Cox et al. \cite{Cox1998}.

\subsection{Isotopisms of algebras}

The concept of isotopism of algebras was introduced by Albert \cite{Albert1942} as a generalization of that of isomorphism. Specifically, two $n$-dimensional algebras $A$ and $A'$ defined over a field $\mathbb{K}$ are {\em isotopic} if there exist three non-singular linear transformations $f$, $g$ and $h$ from $A$ to $A'$ such that
\begin{equation}\label{isot}
f(u)g(v)=h(uv), \text{ for all } u,v\in A.
\end{equation}
The triple $(f,g,h)$ is an {\em isotopism} between the algebras $A$ and $A'$. If $f=g$, then this is called a {\em strong isotopism} and the algebras are said to be {\em strongly isotopic}. If $f=g=h$, then the isotopism constitutes an {\em isomorphism}, which is denoted by $f$ instead of $(f,f,f)$. To be isotopic, strongly isotopic or isomorphic are equivalence relations among algebras. Hereafter, we denote these three relations, respectively, as $\sim$, $\simeq$ and $\cong$.

\vspace{0.1cm}

Let $A$ be an $n$-dimensional algebra over a field $\mathbb{K}$ and let  $\{e_1,\ldots,e_n\}$ be a basis of this algebra. The {\em structure constants} of $A$ are the numbers $c_{ij}^k\in \mathbb{K}$ such that $e_ie_j = \sum_{k=1}^n c_{ij}^k e_k$, for all $1 \leq i, j \leq n$. Its {\em derived algebra} is its subalgebra
$A^2=\{uv\mid\, u,v\in A\}$. The algebra $A$ is {\em abelian} if $A^2$ is trivial, that is, if all its structure constants are zeros. Isotopisms preserve the dimension of derived algebras and hence, that the $n$-dimensional abelian algebra is not isotopic to any other $n$-dimensional algebra. Finally, the {\em left} and {\em right annihilators} of $A$ are respectively defined as the sets
\begin{equation}\label{eq_ann_l}
\mathrm{Ann}_{-}(A)=\{u\in A\mid\, uv=0, \text { for all } v\in A\}.
\end{equation}
\begin{equation}\label{eq_ann_r}
\mathrm{Ann}_{+}(A)=\{u\in A\mid\, vu=0, \text { for all } v\in A\}.
\end{equation}
The intersection of both sets is the {\em annihilator} of $A$,
\begin{equation}\label{eq_ann}
\mathrm{Ann}(A)=\{u\in A\mid\, uv=vu=0, \text { for all } v\in A\}.
\end{equation}

\begin{lemm} \label{lemm_annihilator} Let $(f,g,h)$ be an isotopism between two $n$-dimensional algebras $A$ and $A'$. Then,
\begin{enumerate}[a)]
\item $f(\mathrm{Ann}_{-}(A)) = \mathrm{Ann}_{-}(A')$.
\item $g(\mathrm{Ann}_{+}(A)) = \mathrm{Ann}_{+}(A')$.
\item $f(\mathrm{Ann}_{-}(A))\cap g(\mathrm{Ann}_{+}(A)) = \mathrm{Ann}(A')$.
\end{enumerate}
\end{lemm}

{\bf Proof.} Let us prove (a). Assertion (b) follows similarly and (c) is an immediate consequence of (a) and (b) and the definition of annihilator. Let $u\in A'$ and $v\in f(\mathrm{Ann}_{-}(A))$. Then, $vu=f(f^{-1}(v))g(g^{-1}(u))=h(f^{-1}(v)g^{-1}(u))=h(0)=0$, because $g^{-1}(u)\in A$ and $f^{-1}(v)\in \mathrm{Ann}_{-}(A)$. Hence, $f(\mathrm{Ann}_{-}(A))$ $\subseteq \mathrm{Ann}_{-}(A')$. Now, let $u\in \mathrm{Ann}_{-}(A')$ and $v\in S$. Then, $h(f^{-1}(u)v)=ug(v)=0$. The regularity of $h$ involves that $f^{-1}(u)v=0$. Thus, $u\in f(\mathrm{Ann}_{-}(A))$ and hence, $\mathrm{Ann}_{-}(A')\subseteq f(\mathrm{Ann}_{-}(A))$. \hfill $\Box$

\vspace{0.5cm}

\begin{prop} \label{prop_annihilator} Let $A$ and $A'$ be two $n$-dimensional algebras whose left annihilators coincide, respectively, with their right ones. If both algebras are isotopic, then their annihilators have the same dimension.
\end{prop}

{\bf Proof.} The result is an immediate consequence of Lemma \ref{lemm_annihilator} and the regularity of the components of any isotopism of algebras. \hfill $\Box$

\vspace{0.1cm}

\subsection{Genetic and evolution algebras}

In order to better understand the conceptual meaning of genetic algebras, it is convenient to recall some preliminary concepts in Genetics. A {\em gene} is the molecular unit of hereditary information. This consists of {\em deoxyribonucleic acid} ({\em DNA}), which contains in turn the code to synthesize proteins and determines each one of the attributes that characterize and distinguish each organism. Genes related to a given attribute can have alternative forms, which are called {\em alleles}. Thus, for instance, color of eyes are related to brown, green and blue alleles. Genes are disposed in {\em chromosomes}, which constitute long strands of DNA formed by ordered sequences of genes. The location of alleles related to a given attribute in a chromosome is its {\em locus}, which is preserved by inheritance. Chromosomes carry, therefore, the genetic code of any organism. They also play a main role in the process of reproduction, because the attributes that characterize the offspring are inherited from the alleles that are contained in the chromosomes of the parents. This inheritance depends on the type of organisms under consideration. Thus, for instance, {\em diploid} organisms carry a double set of chromosomes (one of each parent). They reproduce by means of sex cells or {\em gametes}, each of them carrying a single set of chromosomes. The fusion of two gametes of opposite sex gives rise to a {\em zygote}, which contains a double set of chromosomes. Each one of the attributes that characterize the new diploid organism is uniquely determined by the pair of alleles having the same loci in these two chromosomes. There exist distinct laws that regulate, from a probabilistic point of view, the theoretical influence of each one of these two alleles in the final attribute inherited by the offspring. Thus, for instance, the laws of {\em simple Mendelian inheritance} indicate that, for each pair of alleles related to a given attribute, the next generation will inherit with equal frequency both alleles.

\vspace{0.2cm}

Let $\beta=\{e_1,\ldots, e_n\}$ constitute the set of genetically distinct alleles that are related to a given attribute of a population. Then, a {\em genetic algebra} over a field $\mathbb{K}$ that is based on the set $\beta$ is an $n$-dimensional algebra of basis $\beta$ whose structure constants in each product $e_ie_j=\sum_{k=1}^n c_{ij}^ke_k$ refer to the probability that an arbitrary gamete produced by an individual of zygotic type $e_ie_j$ contains the allele $e_k$. Hence, $\sum_{k=1}^n c_{ij}^k=1$, for all $i,j\leq n$. Here, all zygotes have the same fertility and there is an absence of selection. In simple Mendelian inheritance, for instance, we have that $e_ie_j=\frac 12 (e_i + e_j)$, for any pair of alleles $e_i$ and $e_j$. Observe that, if two gametes carry the same allele, then the offspring will inherit it.

\vspace{0.2cm}

In order to deal with asexual reproduction processes, the concept of {\em evolution algebra} was introduced \cite{Tian2008, Tian2006} as a type of $n$-dimensional genetic algebra over a field $\mathbb{K}$ that admits a {\em natural basis} $\{e_1,\ldots, e_n\}$ such that
\begin{itemize}
\item $e_ie_j=0$, if $i\neq j$.
\item $e_ie_i=\sum_{j=1}^n t_{ij}e_j$, for some structure constants $t_{i1},\ldots,t_{in}\in \mathbb{K}$.
\end{itemize}

As a genetic algebra, each basis vector of an evolution algebra constitutes an allele; the product $e_ie_j=0$, for $i\neq j$, represents uniparental inheritance; the product $e_ie_i$ represents self-replication; and each structure constant $t_{ij}$ constitutes the probability that the allele $e_i$ becomes the allele $e_j$ in the next generation. Hereafter, the set of $n$-dimensional evolution algebras over a field $\mathbb{K}$ is denoted as $\mathcal{E}_n(\mathbb{K})$.

\begin{theo}[\cite{Casas2014}]\label{Casas} Every two-dimensional non-abelian complex evolution algebra $A\in\mathcal{E}_2(\mathbb{C})$ is isomorphic to exactly one of the next algebras
\begin{itemize}
\item $\dim A^2=1$:
\begin{itemize}
    \item $E_1:$ $e_1e_1=e_1$.
    \item $E_2:$ $e_1e_1=e_2e_2=e_1$.
    \item $E_3:$ $e_1e_1=e_1+e_2$ and $e_2e_2=-e_1-e_2$.
    \item $E_4:$ $e_1e_1=e_2$.
\end{itemize}
\item $\dim A^2=2$:
\begin{itemize}
    \item $E_{5_{a,b}}:$ $e_1e_1=e_1+ae_2$ and $e_2e_2=be_1+e_2$, where $a,b\in\mathbb{C}$ are such that $ab\neq 1$. Here, $E_{5_{a,b}}\cong E_{5_{b,a}}$, for all $a,b\in\mathbb{C}$.
    \item $E_{6_a}:$ $e_1e_1=e_2$ and $e_2e_2=e_1+ae_2$, where $a\in\mathbb{C}$. If $a,b\in\mathbb{C}\setminus\{0\}$, then $E_{6_a}\cong E_{6_b}$ if and only if $\frac ab=\cos\frac {2k\pi}3+i\sin \frac{2k\pi}3$, for some $k\in\{0,1,2\}$. \hfill $\Box$
\end{itemize}
\end{itemize}
\end{theo}

\vspace{0.1cm}

\subsection{Computational Algebraic Geometry}

Let $I$ be an ideal of a multivariate polynomial ring $\mathbb{K}[X]$. The {\em algebraic set} defined by $I$ is the set $\mathcal{V}(I)$ of common zeros of all its polynomials. The ideal $I$ is {\em zero-dimensional} whenever this set is finite. It is {\em radical} if every polynomial $f\in \mathbb{K}[X]$ belongs to $I$ whenever there exists a natural number $m$ such that $f^m\in I$. The largest monomial of a polynomial in $I$ with respect to a given monomial term ordering is its {\em leading monomial}. The ideal generated by all the leading monomials of $I$ is its {\em initial ideal}. A {\em standard monomial} of $I$ is any monomial that is not contained in its initial ideal. Regardless of the monomial term ordering, if the ideal $I$ is zero-dimensional and radical, then the number of standard monomials in $I$ coincides with the Krull dimension of the quotient ring $\mathbb{K}[X]/I$ and with the number of points of the algebraic set $\mathcal{V}(I)$. This is computed from any Gr\"obner basis of the ideal, that is, from any subset $G$ of polynomials in $I$ whose leading monomials generate its initial ideal. This basis is {\em reduced} if it only has monic polynomials and no monomial of a polynomial in $G$ is generated by the leading monomials of the rest of polynomials in the basis. There exists only one reduced Gr\"obner basis, which can always be computed from Buchberger's algorithm \cite{Buchberger2006}. Similar to the Gaussian elimination on linear systems of equations, this consists of a sequential multivariate division of polynomials. The computation that is required to this end is extremely sensitive to the number of variables.

\vspace{0.1cm}

\begin{theo}[\cite{Gao2009}, Proposition 4.1.1]\label{Gao} Let $\mathbb{F}_q$ be a finite field, with $q$ a power prime. The complexity time that Buchberger's algorithm requires to compute the reduced Gr\"obner bases of an ideal $\langle\, p_1,\ldots, p_m,x_1^q-x_1,\ldots,x^q_n-x_n\,\rangle$ defined over a polynomial ring $\mathbb{F}_q[x_1,\ldots,x_n]$, where $p_1,\ldots,p_m$ are polynomials given in sparse form and have longest length $l$, is $q^{O(n)}+O(m^2l)$. Here, sparsity refers to the number of monomials. \hfill $\Box$
\end{theo}

\vspace{0.5cm}

Computational Algebraic Geometry can be used to determine the isomorphisms and isotopisms between two evolution algebras $A$ and $A'$ in $\mathcal{E}_n(\mathbb{F}_q)$, with respective basis $\{e_1,\ldots,e_n\}$ and $\{e'_1,\ldots,e'_n\}$ and respective structure constants $t_{ij}$ and ${t'}_{ij}$. To this end, we define the sets of variables
$$\mathfrak{F}_n=\{\mathfrak{f}_{ij}\mid\, i,j\leq n\}, \hspace{0.25cm}
\mathfrak{G}_n=\{\mathfrak{g}_{ij}\mid\, i,j\leq n\} \hspace{0.25cm} \text{and} \hspace{0.25cm}
\mathfrak{H}_n=\{\mathfrak{h}_{ij}\mid\, i,j\leq n\},$$
that play the role of the entries in the regular matrices related to a possible isotopism $(f,g,h)$ between the algebras $A$ and $A'$. Here, $\alpha(e_i)=\sum_{j=1}^n \alpha_{ij}e'_j$, for each $\alpha\in\{f,g,h\}$. The next equalities follow in particular from the coefficients of each basis vector $e_l$ in the expression $f(e_i)g(e_j)=h(e_ie_j)$.

\begin{equation}\label{eq_isot_matrices_CAG}
\sum_{k=1}^n \mathfrak{f}_{ik}\mathfrak{g}_{jk}{t'}_{kl} =\begin{cases}
0, \text{ if } i\neq j,\\
\sum_{k=1}^n \mathfrak{h}_{kl}t_{ik}, \text{ otherwise}.
\end{cases}
\end{equation}

\vspace{0.15cm}

\begin{theo}\label{thm_CAG_Isom} Let $A$ and $A'$ be two evolution algebras in $\mathcal{E}_n(\mathbb{F}_q)$, with respective basis $\{e_1,\ldots,e_n\}$ and $\{e'_1,\ldots,e'_n\}$ and respective structure constants $t_{ij}$ and ${t'}_{ij}$. Then,
\begin{enumerate}[a)]
\item The isotopism group between the algebras $A$ and $A'$ is identified with the algebraic set defined by the ideal $I^{\mathrm{Isot}}_{A,A'}$ of $\mathbb{F}_q[\mathfrak{F}_n\cup\mathfrak{G}_n\cup\mathfrak{H}_n]$, which is defined as
$$\langle\, \sum_{k=1}^n \mathfrak{f}_{ik}\mathfrak{g}_{jk}{t'}_{kl}\mid\, i,j,l\leq n; i\neq j\,\rangle + \langle\, \sum_{k=1}^n \mathfrak{f}_{ik}\mathfrak{g}_{ik}{t'}_{kl} - \sum_{k=1}^n \mathfrak{h}_{kl}t_{ik}\mid\, i,l\leq n\,\rangle +$$ $$\langle\,\det(M)^{q-1}-1\mid\, M\in\{F,G,H\}\,\rangle,$$

\noindent where $F$, $G$ and $H$ denote, respectively, the matrices of entries $\{\mathfrak{f}_{ij}\mid i,j\leq n\}$, $\{\mathfrak{g}_{ij}\mid i,j\leq n\}$ and $\{\mathfrak{h}_{ij}\mid i,j\leq n\}$. Besides, $$|\mathcal{V}(I^{\mathrm{Isot}}_{A,A'})|= \mathrm{dim}_{\mathbb{F}_q} (\mathbb{F}_q[\mathfrak{F}_n\cup\mathfrak{G}_n\cup\mathfrak{H}_n]/ I^{\mathrm{Isot}}_{A,A'}).$$

\item The isomorphism group between the algebras $A$ and $A'$ is identified with the algebraic set of the ideal $I^{\mathrm{Isom}}_{A,A'}$ of $\mathbb{F}_q[\mathfrak{F}_n]$, which is defined as
$$\langle\, \sum_{k=1}^n \mathfrak{f}_{ik}\mathfrak{f}_{jk}{t'}_{kl}\mid\, i,j,l\leq n; i\neq j\,\rangle + \langle\, \sum_{k=1}^n \mathfrak{f}^2_{ik}{t'}_{kl} - \sum_{k=1}^n \mathfrak{f}_{kl}t_{ik}\mid\, i,l\leq n\,\rangle +$$
$$\langle\,\det(F)^{q-1}-1\,\rangle,$$
where $F$ denotes the matrix of entries $\{\mathfrak{f}_{ij}\mid\, i,j\leq n\}$. Besides,
$$|\mathcal{V}(I^{\mathrm{Isom}}_{A,A'})|= \mathrm{dim}_{\mathbb{F}_q}(\mathbb{F}_q[\mathfrak{F}_n]/ I^{\mathrm{Isom}}_{A,A'}).$$
\end{enumerate}
\end{theo}

{\bf Proof.} We prove assertion (b), being analogous the reasoning for (a). From (\ref{eq_isot_matrices_CAG}), the generators of the ideal $I^{\mathrm{Isom}}_{A,A'}$ involve each zero of its algebraic set to constitute the entries of the regular matrix of an isomorphism $f$ between the algebras $A$ and $A'$. The result follows from the fact of being this ideal zero-dimensional and radical. Particularly, the ideal $I^{\mathrm{Isom}}_{A,A'}$ is zero-dimensional because its algebraic set is a finite subset of $\mathbb{F}_q^{n^2}$. Besides, from Proposition 2.7 of \cite{Cox1998}, the ideal $I$ is also radical, because, for each $i,j\leq n$, the unique monic generator of $I\cap \mathbb{F}_q[\mathfrak{f}_{ij}]$ is the polynomial $(\mathfrak{f}_{ij})^q-\mathfrak{f}_{ij}$, which is intrinsically included in each ideal of $\mathbb{F}_q[\mathfrak{F}_n]$ and is square-free.  \hfill $\Box$

\vspace{0.5cm}

\begin{coro}\label{coro_CAG_Isom} The complexity times that Buchberger's algorithm requires to compute the reduced Gr\"obner bases of the ideals $I^{\mathrm{Isot}}_{A,A'}$ and $I^{\mathrm{Isom}}_{A,A'}$ in Theorem \ref{thm_CAG_Isom} are, respectively, $q^{O(3n^2)}+O(n^6n!)$ and $q^{O(n^2)}+O(n^6n!)$.
\end{coro}

{\bf Proof.} We prove the result for $I^{\mathrm{Isom}}_{A,A'}$, being analogous the reasoning for $I^{\mathrm{Isot}}_{A,A'}$. The result follows straightforward from Theorem \ref{Gao} once we observe that all the generators of the ideal in Theorem \ref{thm_CAG_Isom} are sparse in $\mathbb{F}_q[\mathfrak{F}_n]$. \hfill $\Box$

\vspace{0.3cm}

Algorithm \ref{alg7} shows how Theorem \ref{thm_CAG_Isom} can be implemented to distribute a subset of $\mathcal{E}_n(\mathbb{F}_q)$ into isotopism and isomorphism classes. Its correctness and termination are based on those of Buchberguer's algorithm \cite{Buchberger2006}.

\begin{algorithm}[h]
\caption{Computation of isomorphism (isotopism, respectively) classes of a set of evolution algebras in $\mathcal{E}_n(\mathbb{F}_q)$.}\label{alg7}
\begin{algorithmic}[1]
\Require A subset $S\subseteq \mathcal{E}_n(\mathbb{F}_q)$.
\Ensure $C$, the set of isomorphism (isotopism, respectively) classes of $S$.
\State $C=\emptyset$.
\While{$S\neq\emptyset$}
\State Take $A\in S$.
\State $S:=S\setminus \{A\}$.
\State $C:=C\cup\{A\}$.
\For {$A' \in S$}
\If{$|\mathcal{V}(I^{\mathrm{Isom}}_{A,A'})|>0$ ($|\mathcal{V}(I^{\mathrm{Isot}}_{A,A'})|>0$, respectively)}
\State $S:=S\setminus \{A'\}$.
\EndIf
\EndFor
\EndWhile
\State \textbf{return} $C$.
\end{algorithmic}
\end{algorithm}

We have implemented Theorem \ref{thm_CAG_Isom} as a procedure called {\em isoAlg} in the open computer algebra system for polynomial computations {\sc Singular} \cite{Decker2016}. This has been included in the library {\em evolution.lib}, which is available online at {\texttt{http://personales.us.es/raufalgan/LS/evolution.lib}}. Having as output the number of isotopisms or that of isomorphisms between two evolution algebras $A$ and $A'$ in $\mathcal{E}_n(\mathbb{F}_q)$, the procedure {\em isoAlg} receives as input the dimension $n$, the order $q$, a pair of lists formed by the structure constants of the algebras under consideration and a positive integer {\em opt} $\leq 2$ that enables us to use the ideal $I^{\mathrm{Isot}}_{A,A'}$, if {\em opt} = 1, or the ideal $I^{\mathrm{Isom}}_{A,A'}$, if {\em opt} = 2. We have made use of this procedure in order to determine all the isotopisms and isomorphisms that appear throughout the paper.

\section{Structure tuples and annihilators of evolution algebras}

This section deals separately with two aspects of evolution algebras that enable us to determine their distribution into isotopism and isomorphism classes: their structure tuples and their annihilators.

\subsection{Structure tuples}

Hereafter, $\{e_1,\ldots,e_n\}$ denotes the natural basis of any evolution algebra in $\mathcal{E}_n(\mathbb{K})$ and $\mathcal{T}_n(\mathbb{K})$ denotes the direct product $\prod_{i=1}^n \langle\,e_1,\ldots,e_n\,\rangle$. The {\em structure tuple} of an evolution algebra in $\mathcal{E}_n(\mathbb{K})$ is the tuple $T=(\mathfrak{t}_1,\ldots,\mathfrak{t}_n)\in \mathcal{T}_n(\mathbb{K})$ where $\mathfrak{t}_i=e_ie_i$, for all $i\leq n$. This algebra is denoted as $A_T$.

\begin{lemm}\label{lemmIsot0_evol} Let $T$ and $T'$ be two structure tuples in $\mathcal{T}_n(\mathbb{K})$ that are equal up to permutation of their components and basis vectors. The evolution algebras $A_T$ and $A_{T'}$ in $\mathcal{E}_n(\mathbb{K})$ are strongly isotopic.
\end{lemm}

{\bf Proof.} Let $T=(\sum_{j=1}^nt_{1j}e_j,\ldots,\sum_{j=1}^nt_{nj}e_j)$ and $T'=(\sum_{j=1}^nt'_{1j}e_j,\ldots,$ $\sum_{j=1}^nt'_{nj}e_j)$. From the hypothesis, there exist two permutations $\alpha$ and $\beta$ of the set $\{1,\ldots,n\}$ such that $t'_{\alpha(i)\beta(j)}=t_{ij}$, for all $i,j<n$. It is then enough to define by linearity the strong isotopism $(f,f,h)$ from $A_T$ to $A_{T'}$ such that $f(e_i)=e_{\alpha(i)}$ and $h(e_i)=e_{\beta(i)}$, for all $i\leq n$. Then, $f(e_i)f(e_j)=0=h(e_ie_j)$, for all $i,j\leq n$ such that $i\neq j$. Besides, for all $i\leq n$, we have that

{\small
$$f(e_i)f(e_i)=e_{\alpha(i)}e_{\alpha(i)} = \sum_{j=1}^n t'_{\alpha(i)j}e_j= \sum_{j=1}^n t'_{\alpha(i)\beta(j)}e_{\beta(j)}= \sum_{j=1}^n t_{ij}e_{\beta(j)}=h(e_ie_i).$$} \hfill $\Box$

\begin{ejemplo}\label{ex_lemmIsot0_evol} From Lemma \ref{lemmIsot0_evol}, we have, for instance, that the evolution algebras $E_1$ and $E_4$ in Theorem \ref{Casas} are strongly isotopic. Specifically, the triple $(\mathrm{Id},\mathrm{Id},h)$, where $h$ switches the basis vectors $e_1$ and $e_2$, is a strong isotopism between both algebras. \hfill $\lhd$
\end{ejemplo}

\vspace{0.2cm}

\begin{prop}\label{propIsot0_evol} Let $T$ be a structure tuple in $\mathcal{T}_n(\mathbb{K})$. There always exists a structure tuple $T'=(\sum_{j=1}^nt'_{1j}e_j,\ldots,\sum_{j=1}^nt'_{nj}e_j)\in \mathcal{T}_n(\mathbb{K})$ such that $A_{T'}$ is strongly isotopic to $A_T$ and the next two conditions hold
\begin{enumerate}[a)]
\item If $t'_{ii}=0$ for some $i\geq 1$, then $t'_{jk}=0$, for all $j,k\geq i$.
\item If $t'_{ii}\neq 0$ for some $i\geq 1$, then $t'_{ij}=0$, for all $j\neq i$.
\end{enumerate}
\end{prop}

{\bf Proof.} Let $T=(\sum_{j=1}^n t_{1j}e_j,\ldots,\sum_{j=1}^nt_{nj}e_j)\in \mathcal{T}_n(\mathbb{K})$. As a first step in the construction of the required structure tuple $T'$, let us consider $T'=T$. From Lemma \ref{lemmIsot0_evol}, any permutation of the components of $T'$, together with any relabeling of the indices of the basis vectors, gives rise to a new evolution algebra in $\mathcal{E}_n(\mathbb{K})$ that is strongly isotopic to $A_T$. Keeping this in mind, we modify conveniently $T'$ so that, if $t'_{ii}=0$, for some $i<n$, then
\begin{itemize}
\item $t'_{ji}=0$, for all $j>i$. Otherwise, we rearrange conveniently from the $i^{th}$ to the $n^{th}$ components of $T'$.
\item $t'_{ij}=0$, for all $j>i$. Otherwise, we permute conveniently the indices of the basis vectors $e_i,\ldots, e_n$.
\end{itemize}

Condition (a) in the statement holds then from the combination of these two assumptions. Now, in order to obtain condition (b), we modify $T'$ so that, for each $i<n$ such that $t'_{ii}\neq 0$, we define by linearity the strong isotopism $(\mathrm{Id},\mathrm{Id},h)$ from $A_{T'}$ in such a way that $h(e_i)=e_i-\frac 1{{t'}_{ii}}(\sum_{j=1}^{i-1} t'_{ij}e_j-\sum_{j=i+1}^n t'_{ij}e_j)$ and $h(e_j)=e_j,$ for all $j\neq i$. Then,
{\small
$$e_ie_i=\mathrm{Id}(e_i)\mathrm{Id}(e_i)=h(e_ie_i)=h(\sum_{j=1}^n t'_{ij}e_j)= \sum_{j=1}^{i-1}t'_{ij}e_j + t'_{ii}h(e_i)+\sum_{j=i+1}^n t'_{ij}e_j=t'_{ii}e_i$$}
and condition (b) holds. \hfill $\Box$

\vspace{0.3cm}

\begin{ejemplo}\label{ex_propIsot0_evol} By following the reasoning exposed in the proof of Proposition \ref{propIsot0_evol}, we obtain that the evolution algebra $E_3$ in Theorem \ref{Casas} is strongly isotopic to the evolution algebra $A_{(e_1,-e_1)}\in\mathcal{E}_2(\mathbb{C})$ by means of the strong isotopism $(\mathrm{Id},\mathrm{Id},h)$, where $h(e_1)=e_1-e_2$ and $h(e_2)=e_2$.

\vspace{0.2cm}

Similarly, any evolution algebra $E_{5_{a,b}}$ in Theorem \ref{Casas} is strongly isotopic to the algebra $E_{5_{0,0}}$. Specifically, if $a$ and $b$ are two complexes numbers such that $ab\neq 1$, then the triple $(\mathrm{Id},\mathrm{Id},h)$ such that $h(e_1)=e_1-ae_2$ and $h(e_2)=e_2$ is a strong isotopism between the evolution algebra $E_{5_{a,b}}$ and the evolution algebra $A_T$ of structure tuple $T=(e_1,be_1+(1-ab)e_2)\in\mathcal{T}_2(\mathbb{C})$. Now, the triple $(\mathrm{Id},\mathrm{Id},h')$ such that $h'(e_2)=\frac 1{1-ab}(e_2-be_1)$ and $h'(e_1)=e_1$ is a strong isotopism between the evolution algebras $A_T$ and $E_{5_{0,0}}$.

\vspace{0.2cm}

Finally, any evolution algebra $E_{6_a}$ in Theorem \ref{Casas} is also strongly isotopic to the algebra $E_{5_{0,0}}$. Specifically, if $a$ is a complex number distinct of zero, then the triple $(f,f,h)$ such that $f$ switches the basis vectors $e_1$ and $e_2$, whereas $h(e_1)=e_1-ae_2$ and $h(e_2)=e_2$, is a strong isotopism between the evolution algebras $E_{6_a}$ and $E_{5_{0,0}}$. \hfill $\lhd$
\end{ejemplo}

\vspace{0.3cm}

We finish this subsection by determining explicitly the distribution of the set $\mathcal{E}_2(\mathbb{C})$ into isotopism classes.

\begin{prop}\label{prop_Ev_Complex} There exist four isotopism classes in $\mathcal{E}_2(\mathbb{C})$. They correspond to the abelian algebra and the evolution algebras $E_1$, $E_2$ and $E_{5_{0,0}}$ in Theorem \ref{Casas}.
\end{prop}

{\bf Proof.} In Examples \ref{ex_lemmIsot0_evol} and \ref{ex_propIsot0_evol} we have seen that $E_1\simeq E_4$, $E_3\simeq A_{(e_1,-e_1)}$ and $E_{5_{a,b}}\simeq E_{6_c}$, for all $a,b,c\in\mathbb{C}$ such that $ab\neq 1$ and $c\neq 0$. Observe now that the triple $(f,f,\mathrm{Id})$ such that $f(e_1)=-ie_2$ and $f(e_2)=e_1$ is a strong isotopism between the algebras $E_2$ and $A_{(e_1,-e_1)}$. Hence, $E_2\simeq E_3$. It is enough to prove, therefore, that the four evolution algebras of the statement are not isotopic. Since the abelian algebra is not isotopic to any other algebra, we can focus on the algebras $E_1$, $E_2$ and $E_{5_{0,0}}$. From Proposition \ref{prop_annihilator}, the former is not isotopic to $E_2$ or $E_{5_{0,0}}$, because $\dim \mathrm{Ann}(E_1)=1\neq 0 = \dim \mathrm{Ann}(E_2)=\dim \mathrm{Ann}(E_{5_{0,0}})$. Finally, since isotopisms preserve the dimension of derived algebras, the evolution algebras $E_2$ and $E_{5_{0,0}}$ are not isotopic, because $E_2^2=\langle\,e_1\,\rangle\subset \langle\,e_1,e_2\,\rangle = E_{5_{0,0}}^2$. \hfill $\Box$

\vspace{0.1cm}

\subsection{Annihilators}

Let $m\leq n$ be a non-negative integer and let $\mathcal{E}_{n;m}(\mathbb{K})$ denote the subset of evolution algebras in $\mathcal{E}_n(\mathbb{K})$ with an $(n-m)$-dimensional annihilator. The set $\{\mathcal{E}_{n;m}(\mathbb{K})\mid\, 0\leq m\leq n\}$ constitutes a partition of the set $\mathcal{E}_n(\mathbb{K})$.

\begin{lemm}\label{lemmEv0a} Every evolution algebra in $\mathcal{E}_{n;m}(\mathbb{K})$ is isomorphic to an evolution algebra in $\mathcal{E}_n(\mathbb{K})$ such that $e_ie_i\neq 0$ if and only if $i\leq m$.
\end{lemm}

{\bf Proof.} Let $A\in \mathcal{E}_{n;m}(\mathbb{K})$. There must exist a subset $S=\{i_1,\ldots,i_m\}\subseteq [n]$ such that $e_ie_i\neq 0$ if and only if $i\in S$. It is then enough to consider the isomorphism that maps, respectively, the basis vectors $e_{i_1},\ldots,e_{i_m}$ to $e_1,\ldots,e_m$ and preserves the rest of basis vectors. \hfill $\Box$

\vspace{0.3cm}

\begin{prop}\label{propEv0} Let $m$ and $m'$ be two distinct non-negative integers less than or equal to $n$. Then, none evolution algebra in $\mathcal{E}_{n;m}(\mathbb{K})$ is isotopic to an evolution algebra in $\mathcal{E}_{n;m'}(\mathbb{K})$.
\end{prop}

{\bf Proof.} The result follows straightforward from Proposition \ref{prop_annihilator} and the fact that $\mathrm{Ann}_{-}(A)=\mathrm{Ann}_{+}(A)=\mathrm{Ann}(A)$, for all $A\in\mathcal{E}_n(\mathbb{K})$. \hfill $\Box$

\vspace{0.3cm}

The next result deals with the distribution of the set $\mathcal{E}_{n;m}(\mathbb{K})$ into isomorphism and isotopism classes, for all positive integer $n\in\mathbb{N}$ and $m\in\{0,1,2\}$. In the statement of the result we make use of the description of the algebras that were exposed in Theorem \ref{Casas} with the exception of dealing here with $n$-dimensional evolution algebras over the field $\mathbb{K}$ instead of two-dimensional complex evolution algebras. Similar abuse of notation is done from here on in order to get a simple and coherent labeling of the evolution algebras that are exposed in the paper.

\begin{prop}\label{propEv1} The next assertions hold.
\begin{enumerate}[a)]
\item The set $\mathcal{E}_{n;0}(\mathbb{K})$ is only formed by the $n$-dimensional abelian algebra.
\item Any evolution algebra in $\mathcal{E}_{1;1}(\mathbb{K})$ is isomorphic to the algebra $E_1$.
\item If $n>1$, then any evolution algebra in $\mathcal{E}_{n;1}(\mathbb{K})$ is isomorphic to the algebra $E_1$ or to the algebra $E_4$.
\item Any evolution algebra in $\mathcal{E}_{n;1}(\mathbb{K})$ is isotopic to the algebra $E_1$.
\item Any evolution algebra in $\mathcal{E}_{2;2}(\mathbb{K})$ is isomorphic to an evolution algebra in $\mathcal{E}_{2;2}(\mathbb{K})$ with natural basis $\{e_1,e_2\}$ such that $e_1e_1\in\{e_1,e_2,e_1+e_2\}$.
\item Any evolution algebra in $\mathcal{E}_{n;2}(\mathbb{K})$ is isotopic to $E_2$ or $E_{5_{0,0}}$.
\end{enumerate}
\end{prop}

{\em Proof.} Let us prove each assertion separately.
\begin{enumerate}[a)]
\item This assertion follows straightforward from the definition of $\mathcal{E}_{n;0}(\mathbb{K})$.
\item Every non-abelian evolution algebra in $\mathcal{E}_1(\mathbb{K})$ is described by a product $e_1e_1=ae_1$, where $a\in\mathbb{K}\setminus\{0\}$. The linear transformation $f$ that maps $e_1$ to $a e_1$ is an isomorphism between this algebra and $E_1$.
\item Let $A$ be an evolution algebra in $\mathcal{E}_{n;1}(\mathbb{K})$ with structure constants $t_{ij}$. From Lemma \ref{lemmEv0a}, we can suppose that $t_{ij}=0$, for all $i>1$. Let $j_0\leq n$ denote the minimum positive integer such that $t_{1j_0}\neq 0$. This exists because $A$ is non-abelian. We can suppose that $j_0\in\{1,2\}$. Otherwise, it is enough to consider the isomorphism that switches the basis vectors $e_2$ and $e_{j_0}$. Let us study both cases.
    \begin{itemize}
        \item If $j_0=1$, then the linear transformation $f$ that is defined such that $f(e_1)=t_{11}e_1-\frac 1{t_{11}}\sum_{j=2}^nt_{1j}e_j$ and $f(e_i)=e_i$, for all $i>1$, is an isomorphism between $A$ and the evolution algebra $E_1$.
        \item If $j_0=2$, then the linear transformation $f$ that is defined such that $f(e_2)=\frac 1{t_{12}}(e_2-\sum_{j=3}^n t_{1j}e_j)$ and $f(e_i)=e_i$, for all $i\neq 2$, is an isomorphism between $A$ and the evolution algebra $E_4$.
    \end{itemize}
\item If $n=1$, then the result follows similarly to (b). Otherwise, it is enough to observe that the triple $(\mathrm{Id},\mathrm{Id},h)$, where $h$ switches the basis vectors $e_1$ and $e_2$, is a strong isotopism between the algebras $E_1$ and $E_4$.
\item Let $A\in\mathcal{E}_{2;2}(\mathbb{K})$. From Lemma \ref{lemmEv0a}, we can suppose the existence of a pair $(a,b)\in\mathbb{K}^2\setminus\{(0,0)\}$ such that $e_1e_1=ae_1+be_2$. If $a\neq 0$, then the algebra $A$ is isomorphic to an evolution algebra such that $e_1e_1=e_1$, if $b=0$, or $e_1e_1=e_1+e_2$, otherwise. To this end, it is enough to consider the isomorphism $f$ such that $f(e_1)=ae_1$ and $f(e_2)=e_2$, if $b=0$, or $f(e_2)=\frac {a^2}b e_2$, otherwise. Further, if $a=0$, then $b\neq 0$ and the algebra $A$ is isomorphic to the evolution algebra such that $e_1e_1=e_2$ by means of the isomorphism that maps $e_2$ to $be_2$ and preserves the basis vector $e_1$.
\item Let $A\in\mathcal{E}_{n;2}(\mathbb{K})$. From Proposition \ref{propIsot0_evol} and Lemma \ref{lemmEv0a}, we can suppose the existence of a pair $(a,b)\in\mathbb{K}^2\setminus\{(0,0)\}$ such that $A$ is strongly isotopic to the evolution algebra in $\mathcal{E}_{n;2}(\mathbb{K})$ with structure tuple $T_1=(ae_1,be_1,0,\ldots,0)$ or $T_2=(ae_1,be_2,0,\ldots,0)$ in $\mathcal{T}_n(\mathbb{K})$. The evolution algebra $A_{T_1}$ is isotopic to $E_2$ by means of the triple $(f,\mathrm{Id},\mathrm{Id})$ such that $f(e_1)=ae_1$, $f(e_2)=be_2$ and $f(e_i)=e_i$, for all $i>2$, whereas the evolution algebra $A_{T_2}$ is strongly isotopic to $E_{5_{0,0}}$ by means of the triple $(\mathrm{Id},\mathrm{Id},h)$ such that $h(e_1)=\frac 1a e_1$, $h(e_2)=\frac 1b e_2$ and $h(e_i)=e_i$, for all $i>2$. \hfill $\Box$
\end{enumerate}

\vspace{0.5cm}

The next theorem, which follows straightforward from Proposition \ref{propEv1}, generalizes Proposition \ref{prop_Ev_Complex} and determines explicitly the distribution of two-dimensional evolution algebras into four isotopism classes, whatever the base field is.

\begin{theo}\label{theo_Ev2} There exist four isotopism classes of two-dimensional evolution algebras over any field. They correspond to the abelian algebra and the evolution algebras $E_1$, $E_2$ and $E_{5_{0,0}}$.\hfill $\Box$
\end{theo}

\vspace{0.5cm}

From the point of view of Genetics, the previous result involves the existence of four distinct classes of asexual diploid organisms up to mutation of their alleles. Their distribution into isomorphism classes requires, however, a further study, which constitutes the final part of this paper.

\section{Isomorphism classes of the set $\mathcal{E}_2(\mathbb{K})$}

As a preliminary study, we focus on the finite field $\mathbb{K}=\mathbb{F}_q$, with $q$ a prime power. Particularly, we have implemented the procedure {\em isoAlg} into Algorithm \ref{alg7}, both of them introduced in Section 2, in order to show in Table \ref{Table_PQ3b_evol} the distribution of the set $\mathcal{E}_2(\mathbb{F}_q)$ into isomorphism classes, for $q\leq 7$. In order to expose the efficiency of our procedure, we also expose in Table \ref{Table_PQ3_rt_evol} the run time and usage memory that are required to compute each distribution. This computation refers to a computer system with an {\em Intel Core i7-2600, with a 3.4 GHz processor and 16 GB of RAM}.

\begin{table}[ht]
\begin{center}
\resizebox{\textwidth}{!}{
\begin{tabular}{lllll}
$q$ & \multicolumn{4}{c}{Structure tuples}\\ \hline
2 & $(0,0)$ & $(e_1,e_1)$ & $(e_2,e_1)$ & $(e_1,e_1+e_2)$\\
\ & $(e_2,0)$ & $(e_1+e_2,e_1+e_2)$ & $(e_2,e_1+e_2)$ & $(e_1,e_2)$\\
\ & $(e_1,0)$ & & & \\
3 & $(0,0)$ & $(e_1+e_2,2e_1+2e_2)$ & $(e_2,e_1+2e_2)$ & $(e_1,e_2)$\\
\ & $(e_2,0)$ & $(e_1,e_1)$ & $(e_1+e_2,2e_1+e_2)$ & \\
\ & $(e_1,0)$ & $(e_2,e_1)$ & $(e_1,e_1+e_2)$ & \\
\ & $(e_1,2e_1)$ & $(e_2,e_1+e_2)$ & $(e_1,2e_1+e_2)$ & \\
5 & $(0,0)$  & $(e_2,e_1)$  & $(e_1+e_2,e_1+3e_2)$  & $(e_1,e_1+e_2)$\\
\ & $(e_2,0)$  & $(e_2,e_1+e_2)$  & $(e_1+e_2,e_1+4e_2)$  & $(e_1,2e_1+e_2)$\\
\ & $(e_1,0)$  & $(e_2,e_1+2e_2)$  & $(e_1+e_2,2e_1+e_2)$  & $(e_1,3e_1+e_2)$\\
\ & $(e_1,e_1)$  & $(e_2,e_1+3e_2)$  & $(e_1+e_2,3e_1+e_2)$  & $(e_1,4e_1+e_2)$\\
\ & $(e_1+e_2,4e_1+4e_2)$  & $(e_2,e_1+4e_2)$  & $(e_1+e_2,2e_1+3e_2)$  & $(e_1,e_2)$\\
\ & $(e_1,2e_1)$  & $(e_1+e_2,e_1+2e_2)$  & $(e_1+e_2,3e_1+2e_2)$  & \\
7 & $(0,0)$ & $(e_2,2e_1+e_2)$ & $(e_1+e_2,e_1+2e_2)$  & $(e_1+e_2,3e_1+5e_2)$\\
\ & $(e_1,0)$  & $(e_2,2e_1+3e_2)$   & $(e_1+e_2,e_1+3e_2)$  & $(e_1+e_2,3e_1+6e_2)$\\
\ & $(e_2,0)$  & $(e_2,3e_1+e_2)$  & $(e_1+e_2,e_1+4e_2)$  & $(e_1+e_2,4e_1+3e_2)$\\
\ & $(e_1,e_1)$  & $(e_2,3e_1+3e_2)$  & $(e_1+e_2,e_1+5e_2)$ & $(e_1+e_2,4e_1+5e_2)$\\
\ & $(e_1,2e_1)$  & $(e_1,e_1+e_2)$  & $(e_1+e_2,e_1+6e_2)$  & $(e_1+e_2,4e_1+6e_2)$\\
\ & $(e_1,3e_1)$  & $(e_1,e_1+2e_2)$  & $(e_1+e_2,2e_1+e_2)$  & $(e_1+e_2,6e_1+3e_2)$\\
\ & $(e_1,e_2)$  & $(e_1,e_1+3e_2)$  & $(e_1+e_2,2e_1+3e_2)$  & $(e_1+e_2,6e_1+5e_2)$\\
\ & $(e_2,e_1)$  & $(e_1,3e_1+e_2)$  & $(e_1+e_2,2e_1+4e_2)$  & $(e_1+e_2,6e_1+6e_2)$\\
\ & $(e_2,e_1+e_2)$  & $(e_1,3e_1+2e_2)$  & $(e_1+e_2,2e_1+5e_2)$  &  \\
\ & $(e_2,e_1+3e_2)$  & $(e_1,3e_1+3e_2)$  & $(e_1+e_2,2e_1+6e_2)$  &  \\
\end{tabular}}
\caption{Distribution into isomorphism classes of the set $\mathcal{E}_2(\mathbb{F}_q)$, for $q\leq 7$.}
\label{Table_PQ3b_evol}
\end{center}
\end{table}

\begin{table}[ht]
\begin{center}
\begin{tabular}{l|cc}
q & Run time & Usage memory \\ \hline
2 & 0 seconds & 0 MB\\
3 & 3 seconds & 0 MB\\
5 & 38 seconds & 80 MB\\
7 & 278 seconds & 1360 MB\\
\end{tabular}
\caption{Run time and memory usage that are required to compute the distribution of the set $\mathcal{E}_2(\mathbb{F}_q)$ into isomorphism classes, for $q\leq 7$.}
\label{Table_PQ3_rt_evol}
\end{center}
\end{table}

The next result follows straightforward from the previous computation.

\begin{theo}\label{theo_EPLR} The sets $\mathcal{E}_2(\mathbb{F}_q)$, with $q\in\{2,3,5,7\}$, are respectively distributed into 9, 13, 23 and 38 isomorphism classes. \hfill $\Box$
\end{theo}

\vspace{0.5cm}

Observe in Table \ref{Table_PQ3b_evol} that the distribution of the set $\mathcal{E}_2(\mathbb{F}_2)$ into nine isomorphism classes agrees with that corresponding to the set $\mathcal{E}_2(\mathbb{C})$ that was exposed in Theorem \ref{Casas}. Nevertheless, this does not hold for finite fields of higher orders. Thus, for instance, the evolution algebra $A_{(e_1,2e_1)}$, which has a one-dimensional derived algebra, is not isomorphic to any of the corresponding four evolution algebras $E_1$ to $E_4$ in $\mathcal{E}_2(\mathbb{F}_q)$, for $q>2$. A further study that generalizes the result of Casas et al. \cite{Casas2014} is then required for a general base field $\mathbb{K}$. From Proposition \ref{propEv1}, we can focus on the distribution of the set $\mathcal{E}_{2,2}(\mathbb{K})$ into isomorphism classes and, more specifically, on those two-dimensional evolution algebras with natural basis $\{e_1,e_2\}$ such that $e_1e_1\in\{e_1,e_2,e_1+e_2\}$. To this end, we propose here the use of Computational Algebraic Geometry. Particularly, we eliminate in Theorem \ref{thm_CAG_Isom} those generators of the ideal $I^{\mathrm{Isom}}_{A,A'}$ that are referred to the determinants of the matrices $F$, $G$ and $H$. This reduces the corresponding complexity time that is exposed in Corollary \ref{coro_CAG_Isom} to $q^{O(n^2)} + O(n^8)$ and gives enough information to analyze a case study on which base the possible isomorphisms between two given evolution algebras.

\vspace{0.2cm}

Let $A=A_{(ae_1+be_2,ce_1+de_2)}$ and $A'=A_{(\alpha e'_1+\beta e'_2,\gamma e'_1+\delta e'_2)}$ be two isomorphic evolution algebras in $\mathcal{E}_{2,2}(\mathbb{K})$ with respective natural bases $\{e_1,e_2\}$ and $\{e'_1,e'_2\}$. Let $f$ be an isomorphism between both algebras with a related non-singular matrix $F=(f_{ij})$ such that $f(e_i)=f_{i1}e'_1+f_{i2}e'_2$, for all $i\in\{1,2\}$. The implementation of the procedure {\em isoAlg} enables us to ensure that, whatever the base field is, the reduced Gr\"obner basis of the ideal in Theorem \ref{thm_CAG_Isom} related to the isomorphism group between the evolution algebras $A$ and $A'$ involves in particular that
\begin{equation}\label{eq_ev_1}
\begin{cases}
(ad - bc)f_{11}f_{21}=0,\\
(ad - bc)f_{12}f_{22}=0.
\end{cases}
\end{equation}

From the previous conditions, we can distinguish two cases depending on the fact of being $ad=bc$ or $ad\neq bc$. They refer, respectively, to two-dimensional evolution algebras with a one- or two-dimensional derived algebra. Recall in this regard that any isomorphism between two algebras preserves the dimension of their corresponding derived algebras. In the next two subsections we study separately each one of the two mentioned cases.

\subsection{One-dimensional derived algebra ($ad=bc$)}

In this subsection, $f$ is an isomorphism of regular matrix $F=(f_{ij})$ between the algebras $A=A_{(ae_1 + be_2,c e_1 + d e_2)}$ and $A'=A_{(\alpha e'_1+\beta e'_2,\gamma e'_1 + \delta e'_2)}$ in $\mathcal{E}_{2,2}(\mathbb{K})$ such that $ad=bc$ and $\alpha\delta=\beta\gamma$. From assertion (e) in Proposition \ref{propEv1}, we can suppose that $a,b,\alpha,\beta\in\{0,1\}$. Firstly, suppose $A=A_{(e_1,c e_1)}$, with $c\in\mathbb{K}\setminus\{0\}$. Assertion (e) in Proposition \ref{propEv1} gives rise to the next case study

\begin{itemize}
\item {\bf Case 1}. $A'=A_{(e'_1,\gamma e'_1)}$, with $\gamma\in\mathbb{K}\setminus\{0\}$.

    The identification of coefficients of a same basis vector in each one of the equalities $f(e_ie_j)=f(e_i)f(e_j)$, for all $i,j\leq 2$, involves that $f$ is an isomorphism between the two algebras under consideration if and only if $f_{11}f_{21}=f_{12}f_{22}=0$. The regularity of the matrix $F$ involves that $f_{11}=f_{22}=0$ or $f_{21}=f_{12}=0$. In the first case, we obtain that $f_{21}$ must be zero, what is a contradiction with the regularity of the matrix $F$. In the second case, we obtain that $c=\gamma f_{22}^2$. This fact enables us to ensure that $A_{(e_1,c e_1)}\cong A_{(e_1,c m^2 e_1)}$, for all $c,m\in\mathbb{K}\setminus\{0\}$.

\item {\bf Case 2}. $A'=A_{(e'_2,\delta e'_2)}$, with $\delta\in\mathbb{K}\setminus\{0\}$.

    The computation of the corresponding reduced Gr\"obner basis, which has previously been mentioned, related to these assumptions enables us to ensure that
    \begin{equation}\label{eq_ev2_c12}
    \begin{cases}
    f_{11}=f_{22}=0,\\
    f_{12} = 1 / \delta,\\
    f_{21}^2 = c / \delta.
    \end{cases}
    \end{equation}
    If we take $f_{21}=1$, then we can ensure in particular that $A_{(e_2,c e_2)}\cong A_{(e_1,c e_1)}$, for all $c\in\mathbb{K}\setminus\{0\}$.

\item {\bf Case 3}. $A'=A_{(e'_1 + e'_2,\gamma (e'_1 + e'_2))}$, with $\gamma\in\mathbb{K}\setminus\{0\}$.

    From the reduced Gr\"obner basis related to this case, we deduce that
    \begin{equation}\label{eq_ev2_c13}
    \begin{cases}
    \gamma\neq -1,\\
    f_{11} = f_{12} = 1/(\gamma + 1),\\
    f_{21} = -\gamma f_{22},\\
    c =\gamma (\gamma + 1)^2 f_{22}^2.
    \end{cases}
    \end{equation}
    Particularly, the determinant of the matrix $F$ coincides with $f_{22}$, which must be distinct of zero. As a consequence, the algebras $A_{(e_1 + e_2,c (e_1 + e_2))}$ and $A_{(e_1,c (c + 1)^2 e_1)}$ are isomorphic, for all $c\in\mathbb{K}\setminus\{0,-1\}$.
\end{itemize}

From the previous case study, the case $A=A_{(e_2,d e_2)}$, with $d\in\mathbb{K}\setminus\{0\}$, can be referred to Case 1.2, because $A_{(e_2,d e_2)}\cong A_{(e_1,d e_1)}$. Besides, the case $A=A_{(e_1+e_2,c (e_1+e_2))}$, with $c\in\mathbb{K}\setminus\{0\}$, can be referred to Case 1.3 except for the case $c=-1$, that is, except for the evolution algebra $A_{(e_1+e_2,-(e_1+e_2))}$. The next results gather together what we have just exposed in the previous case study.

\begin{prop}\label{prop_clev1} The next assertions hold in the set $\mathcal{E}_{2,2}(\mathbb{K})$.
\begin{enumerate}[a)]
\item  $A_{(e_1,c e_1)}\cong A_{(e_1,c m^2 e_1)}$ for all $c,m\in\mathbb{K}\setminus\{0\}$.
\item $A_{(e_2,c e_2)}\cong A_{(e_1,c e_1)}$, for all $c\in\mathbb{K}\setminus\{0\}$.
\item $A_{(e_1 + e_2,c (e_1 + e_2))}\cong A_{(e_1,c (c + 1)^2 e_1)}$, for all $c\in\mathbb{K}\setminus\{0,-1\}$. \hfill $\Box$
\end{enumerate}
\end{prop}

\vspace{0.5cm}

\begin{theo}\label{theo_clev1} Any two-dimensional evolution algebra in $\mathcal{E}_{2,2}(\mathbb{K})$ with a one-dimensional derived algebra is isomorphic to exactly one of the next algebras
\begin{itemize}
\item $A_{(e_1,c e_1)}$, with $c\in\mathbb{K}\setminus\{0\}$. Here, $A_{(e_1,c e_1)}\cong A_{(e_1,\gamma e_1)}$ if and only if $\gamma = c m^2$ for some $m\in\mathbb{K}\setminus\{0\}$.
\item $A_{(e_1 + e_2,-e_1 - e_2)}$. \hfill $\Box$
\end{itemize}
\end{theo}

\subsection{Two-dimensional derived algebra ($ad\neq bc$)}

Let us focus now on the case in which the evolution algebra under consideration, $A=A_{(ae_1 + be_2,c e_1 + d e_2)}\in\mathcal{E}_{2,2}(\mathbb{K})$, is such that $ad\neq bc$.
\begin{lemm}\label{lemm_clev0} Let $A=A_{(ae_1+be_2,ce_1+de_2)}$ be a two-dimensional evolution algebra in $\mathcal{E}_{2,2}(\mathbb{K})$ such that $ad\neq bc$. Then, any isomorphism from $A$ to another evolution algebra in $\mathcal{E}_{2,2}(\mathbb{K})$, with related regular matrix $F=(f_{ij})$, holds that $f_{11}=f_{22}=0$ or $f_{12}=f_{21}=0$.
\end{lemm}

{\bf Proof.} The result follows straightforward from both conditions in (\ref{eq_ev_1}) and the regularity of the matrix $F$. \hfill $\Box$

\vspace{0.5cm}

Let us study each case in Lemma \ref{lemm_clev0} separately. Here, $f$ is an isomorphism of regular matrix $F=(f_{ij})$ between a pair of evolution algebras $A=A_{(ae_1 + be_2,c e_1 + d e_2)}$ and $A'=A_{(\alpha e'_1+\beta e'_2,\gamma e'_1 + \delta e'_2)}$ in $\mathcal{E}_{2,2}(\mathbb{K})$, where $ad\neq bc$ and $\alpha\delta\neq \beta\gamma$.

\begin{enumerate}
\item {\bf Case 1}. $f_{11}=f_{22}=0$.

Here, $f(e_1)=f_{12}e'_2$ and $f(e_2)=f_{21}e'_1$, where $f_{12}\neq 0\neq f_{21}$. Similarly to the reasoning exposed in the case study of the previous subsection, the identification of coefficients of a same basis vector in the equalities $f(e_ie_j)=f(e_i)f(e_j)$, for all $i,j\leq 2$, involves that $f$ is an isomorphism between the two algebras under consideration if and only if
    \begin{equation}\label{eq_ev_c1}
    \begin{cases}
    a = \delta f_{12},\\
    b f_{21}= \gamma f_{12}^2,\\
    c f_{12}=\beta f_{21}^2,\\
    d = \alpha f_{21}.
    \end{cases}
    \end{equation}
The regularity of the matrix $F$ implies that a coefficient $a$, $b$, $c$ or $d$ is zero in the structure tuple of the algebra $A$ if and only if the respective coefficient $\delta$, $\gamma$, $\beta$ or $\alpha$ is zero in the structure tuple of $A'$. Now, assertion (e) in Proposition \ref{propEv1} enables us to focus on the following cases for the evolution algebras $A$ and $A'$ under the conditions of Lemma \ref{lemm_clev0}.
    \begin{itemize}
    \item {\bf Case 1.1}. $A=A_{(e_1,c e_1 + d e_2)}$ and $A'=A_{(e'_1, \delta e'_2)}$, where $d\neq 0\neq \delta$. From (\ref{eq_ev_c1}), we have that
        \begin{equation}\label{eq_ev_c11}
        \begin{cases}
        f_{12} = 1 / \delta ,\\
        f_{21} = d,\\
        c=0.
        \end{cases}
        \end{equation}
        Hence, $A_{(e_1,d e_2)}\cong A_{(e_1,e_2)}$, for all $d\in\mathbb{K}\setminus\{0\}$.

    \item {\bf Case 1.2}. $A=A_{(e_1,c e_1 + d e_2)}$ and $A'=A_{(e'_1 + e'_2, \delta e'_2)}$, where $d\neq 0\neq \delta$. From (\ref{eq_ev_c1}), we have that
        \begin{equation}\label{eq_ev_c112}
        \begin{cases}
        f_{12} = 1 / \delta ,\\
        f_{21} = d,\\
        \delta = c / d^2.
        \end{cases}
        \end{equation}
        Hence, $A_{(e_1 + e_2, d e_2)}\cong A_{(e_1,d e_1 + e_2)}$, for all $d\in\mathbb{K}\setminus\{0\}$.

    \item {\bf Case 1.3}. $A=A_{(e_2,c e_1 + d e_2)}$  and $A'=A_{(e'_2,\gamma e'_1)}$, where $c\neq 0\neq \gamma$. From (\ref{eq_ev_c1}), we have that
        \begin{equation}\label{eq_ev_c13}
        \begin{cases}
        f_{21}= \gamma f_{12}^2,\\
        c=\gamma^2f_{12}^3,\\
        d = 0.
        \end{cases}
        \end{equation}
        Hence, $A_{(e_2, c^2m^3 e_1)}\cong A_{(e_2,c e_1)}$, for all $c,m\in\mathbb{K}\setminus\{0\}$.

    \item {\bf Case 1.4}. $A=A_{(e_2,c e_1 + d e_2)}$ and $A'=A_{(e'_1 + e'_2,\gamma e'_1)}$, where $c$, $d$ and $\gamma$ are all of them distinct of zero. From (\ref{eq_ev_c1}), we have that
        \begin{equation}\label{eq_ev_c14}
        \begin{cases}
        f_{12}=d^2/c,\\
        f_{21} = d,\\
        c^2= \gamma d^3.
        \end{cases}
        \end{equation}
        Hence, $A_{(e_1 + e_2,c e_1)}\cong A_{(e_2,\frac 1c (e_1 + e_2))}$, for all $c\in\mathbb{K}\setminus\{0\}$.

    \item {\bf Case 1.5}. $A=A_{(e_1+e_2,c e_1 + d e_2)}$ and $A'=A_{(e'_1 + e'_2,\gamma e'_1 + \delta e'_2)}$, where $c\neq d$ and $c$, $d$, $\gamma$ and $\delta$ are all of them distinct of zero. From (\ref{eq_ev_c1}), we have that
        \begin{equation}\label{eq_ev_c15}
        \begin{cases}
        f_{12} = 1/\delta,\\
        f_{21}= d,\\
        \gamma = c^2/d^3,\\
        \delta = c/d^2.
        \end{cases}
        \end{equation}
        Hence, $A_{(e_1 + e_2,c e_1 + d e_2)}\cong A_{(e_1 + e_2,\frac c{d^2} (\frac cd e_1 + e_2))}$, for all $c,d\in\mathbb{K}\setminus\{0\}$ such that $c\neq d$.
    \end{itemize}

\item {\bf Case 2}. $f_{12}=f_{21}=0$.

Here, $f(e_1)=f_{11}e'_1$ and $f(e_2)=f_{22}e'_2$, where $f_{11}\neq 0\neq f_{22}$. Similarly to the previous case, the identification of coefficients of a same basis vector in the equalities $f(e_ie_j)=f(e_i)f(e_j)$, for all $i,j\leq 2$, involves that
    \begin{equation}\label{eq_ev_c2}
    \begin{cases}
    a = \alpha f_{11},\\
    f_{22}b=\beta f_{11}^2,\\
    f_{11}c=\gamma f_{22}^2,\\
    d = \delta f_{22}.
    \end{cases}
    \end{equation}
    Again from the regularity of the matrix $F$, we have that $a$, $b$, $c$ or $d$ is zero if and only if $\alpha$, $\beta$, $\gamma$ or $\delta$ is zero, respectively. From assertion (e) in Proposition \ref{propEv1}, we consider the following case study.
    \begin{itemize}
    \item {\bf Case 2.1}. $A=A_{(e_1,c e_1 + d e_2)}$ and $A'=A_{(e'_1,\gamma e'_1 + \delta e'_2)}$, where $d\neq 0\neq \delta$. From (\ref{eq_ev_c2}), we have that
    \begin{equation}\label{eq_ev_c21}
    \begin{cases}
    f_{11} = 1,\\
    f_{22} = d / \delta\\
    c\delta^2 = d^2\gamma.
    \end{cases}
    \end{equation}
    Hence, $A_{(e_1,c e_1 + d e_2)}\cong A_{(e_1,\gamma e_1 + \delta e_2)}$ for all $c,d,\gamma,\delta\in\mathbb{K}$ such that $d\neq 0\neq \delta$ and $c\delta^2 = d^2\gamma$.

    \item {\bf Case 2.2}. $A=A_{(e_2,c e_1 + d e_2)}$  and $A'=A_{(e'_2,\gamma e'_1 + \delta e'_2)}$, where $c\neq 0\neq \gamma$. From (\ref{eq_ev_c2}), we have that
    \begin{equation}\label{eq_ev_c22}
    \begin{cases}
    f_{11}^3 = c / \gamma,\\
    f_{22}=f_{11}^2,\\
    d = \delta f_{22}.
    \end{cases}
    \end{equation}
    Hence, $A_{(e_2,c e_1 + d e_2)}\cong A_{(e_2,\frac c{m^3} e_1 + \frac d{m^2} e_2)}$ for all $c,d,m\in\mathbb{K}\setminus\{0\}$.

    \item {\bf Case 2.3}. $A=A_{(e_1+e_2,c e_1 + d e_2)}$ and $A'=A_{(e'_1+e'_2,\gamma e'_1 + \delta e'_2)}$, where $c\neq d$. From (\ref{eq_ev_c2}), we have that $f$ is the trivial isomorphism and that $A'$ coincides with $A$.
    \end{itemize}
\end{enumerate}

\vspace{0.3cm}

The next results gather together what we have just exposed in the previous case study.

\begin{prop}\label{prop_clev0} The next assertions hold in the set $\mathcal{E}_{2,2}(\mathbb{K})$.
\begin{enumerate}[a)]
\item $A_{(e_1,d e_2)}\cong A_{(e_1,e_2)}$, for all $d\in\mathbb{K}\setminus\{0\}$.
\item $A_{(e_1 + e_2, d e_2)}\cong A_{(e_1,d e_1 + e_2)}$, for all $d\in\mathbb{K}\setminus\{0\}$.
\item $A_{(e_2, c^2m^3 e_1)}\cong A_{(e_2,c e_1)}$, for all $c,m\in\mathbb{K}\setminus\{0\}$.
\item $A_{(e_1 + e_2,c e_1)}\cong A_{(e_2,\frac 1c (e_1 + e_2))}$, for all $c\in\mathbb{K}\setminus\{0\}$.
\item $A_{(e_1 + e_2,c e_1 + d e_2)}\cong A_{(e_1 + e_2,\frac c{d^2} (\frac cd e_1 + e_2))}$, for all $c,d\in\mathbb{K}\setminus\{0\}$ such that $c\neq d$.
\item $A_{(e_1,c e_1 + d e_2)}\cong A_{(e_1,\gamma e_1 + \delta e_2)}$ for all $c,d,\gamma,\delta\in\mathbb{K}$ such that $d\neq 0\neq \delta$ and $c\delta^2 = d^2\gamma$.
\item $A_{(e_2,c e_1 + d e_2)}\cong A_{(e_2,\frac c{m^3} e_1 + \frac d{m^2} e_2)}$ for all $c,d,m\in\mathbb{K}\setminus\{0\}$. \hfill $\Box$
\end{enumerate}
\end{prop}

\vspace{0.5cm}

\begin{theo}\label{theo_clev0} Any two-dimensional evolution algebra in $\mathcal{E}_{2,2}(\mathbb{K})$ with a two-dimensional derived algebra is isomorphic to exactly one of the next algebras
\begin{itemize}
\item $A_{(e_1, c e_1 + d e_2)}$, with $d\neq 0$.

Here, $A_{(e_1,c e_1 + d e_2)}\cong A_{(e_1,\gamma e_1 + \delta e_2)}$ if and only if $c\delta^2 = d^2\gamma$.
\item $A_{(e_2, c e_1 + d e_2)}$, with $c\neq 0$.

Here, $A_{(e_2,c e_1 + d e_2)}\cong A_{(e_2,\gamma e_1 + \delta e_2)}$ if and only if  there exists an element $m\in\mathbb{K}\setminus\{0\}$ such that $c = \gamma m^3$ and $d=\delta m^2$, or, $c=\gamma^2 m^3$ and $d=\delta=0$.
\item $A_{(e_1 + e_2, c e_1 + d e_2)}$, with $c\neq 0\neq d$.

Here, $A_{(e_1 + e_2,c e_1 + d e_2)}\cong A_{(e_1 + e_2,\gamma e_1 + \delta e_2)}$ if and only if $\gamma=c^2/d^3$ and $\delta=c/d^2$. \hfill $\Box$
\end{itemize}
\end{theo}

\vspace{0.5cm}

\section{Conclusion and further studies}

This paper has dealt with the set $\mathcal{E}_n(\mathbb{K})$ of $n$-dimensional evolution algebras over a field $\mathbb{K}$, whose distribution into isotopism classes is uniquely related with mutations in non-Mendelian Genetics. Particularly, we have focused on the two-dimensional case, which is related to the asexual reproduction processes of diploid organisms. We have proved that the set $\mathcal{E}_2(\mathbb{K})$ is distributed into four isotopism classes, whatever the base field is, and we have characterized its isomorphism classes. Similar case studies to those ones that have been exposed in this paper are established as further work in order to deal with the distribution into isotopism and isomorphism classes of evolution algebras of dimension $n>2$ over any base field.

\end{document}